\documentclass[11pt]{article}
\usepackage[dvips]{graphicx}
\usepackage{amsfonts}
\usepackage[english]{babel}
\usepackage[dvips]{color}
\usepackage{amsmath}
\usepackage{amsthm}
\usepackage{amssymb}
\usepackage{amsfonts}
\usepackage{xcolor}
\usepackage{caption}
\usepackage{amsthm}
\usepackage{xcolor}

\def\vsn{\par \noindent}
\def\vvsn{\vskip 0.5truecm\noindent}
\begin{document}
\begin{center}
{\Large {\bf
{Some fractional integral and derivatives revised}}}
\vvsn
\\ {{\bf Juan Luis Gonzales-Santander}$^a$, 
{\bf Francesco Mainardi}$^b$
}
\vvsn
 $^a$
 Orcid:{0000-0001-5348-4967}
 \\ Department  of Mathematics, University of Oviedo,\\
C Leopoldo Calvo Sotelo 18, E-33007 Oviedo, Spain\\
\vsn
$^b$ 
Orcid:{0000-0003-4858-7309}
\\ Department of Physics $\&$ Astronomy,
University of Bologna, and INFN,
\\ Via Irnerio 46, I-40126 Bologna, Italy
\end{center}
\begin{abstract}
In the most common literature about fractional calculus, we find that $%
_{a}D_{t}^{\alpha }f\left( t\right) =\,_{a}I_{t}^{-\alpha }f\left( t\right) $
is assumed implicitly in the tables of fractional integrals and derivatives.
However, this is not straightforward from the definitions of $%
_{a}I_{t}^{\alpha }f\left( t\right) $ and $_{a}D_{t}^{\alpha }f\left(
t\right) $. In this sense, we prove that $_{0}D_{t}^{\alpha }f\left(
t\right) =\,_{0}I_{t}^{-\alpha }f\left( t\right) $ is true for $f\left(
t\right) =t^{\nu -1}\log t$, and $f\left( t\right) =e^{\lambda t}$, despite the fact that these
derivations are highly non-trivial. Moreover, the corresponding formulas for
$_{-\infty }D_{t}^{\alpha }\left\vert t\right\vert ^{-\delta }$ and $%
_{-\infty }I_{t}^{\alpha }\left\vert t\right\vert ^{-\delta }$ found in the
literature are incorrect; thus, we derive the correct ones, proving in turn
that $_{-\infty }D_{t}^{\alpha }\left\vert t\right\vert ^{-\delta
}=\,_{-\infty }I_{t}^{-\alpha }\left\vert t\right\vert ^{-\delta }$ holds
true
\vspace{3mm}
\\
{\bf Publication:} Mathematics (MDPI) 2024, No12, 2786/1--13. 
\\  DOI:  10.3390/math12172786 (Open access)
\end{abstract}



\section{Introduction}
From the very beginning of the invention of differential calculus, important
inquiries regarding the significance of the non-integer operations of integrals
and derivatives calculus were brought up. In~this sense, it is well known
that Leibniz initially introduced a symbolic approach and employed the
notation $d^{n}y/dt^{n}=D^{n}y$ to represent the $n$-th derivative, with~$n$
being a non-negative integer. However, L'Hospital asked in a letter to
Leibniz dated in $1695$~\cite{Leibniz}: ``What if $n$ is $1/2$?'' Leibniz
replied, ''It will lead to a paradox.'' But~he added, ``From this apparent
paradox, one day useful consequences will be drawn''. From~this initial
``paradox'', fractional calculus was developed through contributions from
mathematicians such as Euler, Lagrange, Laplace, Fourier, and~others during
the $18$th and early $19$th centuries (a comprehensive summary of its
historical progression can be found in \cite{Miller} \mbox{[Chap. I]}
).
Despite the efforts of these great mathematicians, a~satisfactory expression for the
generalization of integration to fractional powers was not developed until
the mid-$19$th century, through the work of Liouville~\cite{Liouville}.
However, it is worth noting that Abel set the notation that was used later by Liouville (and also used nowadays) for fractional-order integration when solving the generalization of the tautochrone problem (see~\cite{Abel} and the references therein).
For a rigorous understanding of fractional calculus as
a theory involving operators of integration and differentiation of arbitrary
order, we recommend the book 
by Samko, Kilbas, and~Marichev~\cite{BibleFC}.
\vsn
{\bf Definition 1.}
 {Riemann--Liouville fractional integral}
 \\
For $\alpha >0$ \cite{Ederly} ([Chap. XIII])%
\begin{equation}
_{a}I_{t}^{\alpha }f\left( t\right) =\frac{1}{\Gamma \left( \alpha \right) }%
\int_{a}^{t}\left( t-\tau \right) ^{\alpha -1}f\left( \tau \right) \,d\tau .
\label{Definition_RL}
\end{equation}
\vsn
{\bf Remark 1.}
Usually, $a=0$ in many textbooks and applications \cite{Mainardi} [Eqn. 1.2].
However, when $a\rightarrow -\infty $ in (\ref{Definition_RL}), we obtain
the Weyl fractional integral, i.e.,
\begin{equation}
_{-\infty }I_{t}^{\alpha }f\left( t\right) =\frac{1}{\Gamma \left( \alpha
\right) }\int_{-\infty }^{t}\left( t-\tau \right) ^{\alpha -1}f\left( \tau
\right) \,d\tau .  \label{Definition_Weyl_FI}
\end{equation}
\vsn
It is worth noting that in \cite{Ederly} [Chap. XIII], we find other
definition of the Weyl fractional integral, i.e.,
\begin{equation*}
\mathfrak{B}_{\alpha }\left\{ f\left( \tau \right) ;t\right\} =\frac{1}{%
\Gamma \left( \alpha \right) }\int_{t}^{\infty }\left( \tau -t\right)
^{\alpha -1}f\left( \tau \right) d\tau .
\end{equation*}%
It is easy 
 to prove that
\begin{equation*}
_{-\infty }I_{t}^{\alpha }f\left( t\right) =\mathfrak{B}_{\alpha }\left\{
f\left( -\tau \right) ;-t\right\} .
\end{equation*}
\vsn
As pointed out in \cite{Mainardi} [Sect. 1.2], one is tempted to substitute $%
\alpha $ with $-\alpha $ in (\ref{Definition_RL}) in order to obtain a
definition for the fractional derivative $_{a}D_{t}^{\alpha }f(t)$.
Nevertheless, some care is required in the integration for this
generalization, and~the theory of generalized functions has to be invoked.
In order to avoid the use of generalized functions, we find in \cite{Polubny} \mbox{[Sect.
2.3.3]} the following definition:
\vsn
{\bf Definition 2.}Riemann--Liouville fractional derivative.
\\
For $m\in
\mathbb{N}
$, and~$\alpha >0$
\begin{equation}
_{a}D_{t}^{\alpha }f\left( t\right) =\left\{
\begin{array}{ll}
\displaystyle%
\frac{1}{\Gamma \left( m-\alpha \right) }\frac{d^{m}}{dt^{m}}\int_{a}^{t}%
\frac{f\left( \tau \right) }{\left( t-\tau \right) ^{\alpha +1-m}}d\tau , &
m-1<\alpha <m, \\
\displaystyle%
\frac{d^{m}}{dt^{m}}f\left( t\right) , & \alpha =m.%
\end{array}%
\right.  \label{Derivative_RL_def}
\end{equation}
\vsn
{\bf Remark 2.}
Usually, $a=0$ in many textbooks and applications \cite{Mainardi} [Eqn. 1.13b]%
. However, when $a\rightarrow -\infty $ in (\ref{Derivative_RL_def}), we
obtain the Weyl fractional derivative \cite{Mainardi} [Eqn. 1.108], i.e.,~for $%
m\in
\mathbb{N}
$, and~$\alpha >0$
\begin{equation}
_{-\infty }D_{t}^{\alpha }f\left( t\right) =\left\{
\begin{array}{ll}
\displaystyle%
\frac{1}{\Gamma \left( m-\alpha \right) }\frac{d^{m}}{dt^{m}}\int_{-\infty
}^{t}\frac{f\left( \tau \right) }{\left( t-\tau \right) ^{\alpha +1-m}}d\tau
, & m-1<\alpha <m, \\
\displaystyle%
\frac{d^{m}}{dt^{m}}f\left( t\right) , & \alpha =m.%
\end{array}%
\right.  \label{Derivative_Weyl_def}
\end{equation}
\vsn
Nonetheless, in~the existing literature, we find tables of Riemann--Liouville
fractional derivatives, wherein they just substitute $\alpha $ with $-\alpha $
in the corresponding Riemann--Liouville fractional integral. For~instance, if~$E_{\mu ,\nu }\left( z\right) $ denotes the Mittag--Leffler function, and~$%
\psi \left( z\right) $ denotes the digamma function, we find in \cite{Ederly} [Eqn.
13.1.(24)] and \cite{Magin} [Table IV.1], respectively,%
\begin{eqnarray}
_{0}I_{t}^{\alpha }\left[ t^{\nu -1}\log t\right] &=&\frac{t^{\alpha +\nu
-1}\,\,\Gamma \left( \nu \right) }{\Gamma \left( \alpha +\nu \right) }\left[
\log t+\psi \left( \nu \right) -\psi \left( \alpha +\nu \right) \right] ,
\label{Fractional_integral_log} \\
_{0}I_{t}^{\alpha }e^{\lambda t} &=&t^{\alpha }E_{1,1+\alpha }\left( \lambda
t\right) ,  \label{Fractional_integral_exp}
\end{eqnarray}%
Nevertheless, in~\cite{Polubny} [Appendix], we find that $\alpha $ is
changed by $-\alpha $ in (\ref{Fractional_integral_log}) and (\ref%
{Fractional_integral_exp}) in order to obtain the corresponding fractional
derivatives, i.e.,~$_{0}D_{t}^{\alpha }\left[ t^{\nu -1}\log t\right] $ and $%
_{0}D_{t}^{\alpha }e^{\lambda t}$.

Also, we find in \cite{Mainardi} [Eqn. 1.112]%
\begin{equation}
\,_{-\infty }D_{t}^{\alpha }\left\vert t\right\vert ^{-\delta }=\frac{\Gamma
\left( \delta +\alpha \right) }{\Gamma \left( \delta \right) }\,\left\vert
t\right\vert ^{-\alpha -\delta }=\,_{-\infty }I_{t}^{-\alpha }\left\vert
t\right\vert ^{-\delta }.  \label{Mainardi_power}
\end{equation}%
However, according to our numerical experiments, it seems that (\ref%
{Mainardi_power}) does not hold true. Consequently, the~aim of this paper is
twofold. On~the one hand, we want to justify that
\begin{eqnarray*}
_{0}D_{t}^{\alpha }\left[ t^{\nu -1}\log t\right] &=&\,_{0}I_{t}^{-\alpha }%
\left[ t^{\nu -1}\log t\right] , \\
_{0}D_{t}^{\alpha }e^{\lambda t} &=&\,_{0}I_{t}^{-\alpha }e^{\lambda t},
\end{eqnarray*}%
from the definitions given in (\ref{Definition_RL}) and (\ref%
{Derivative_RL_def}). We will see that these proofs are highly non-trivial.
On the other hand, we would like to calculate the Weyl fractional integral
and derivative for the power function, i.e.,~$_{-\infty }I_{t}^{\alpha
}\left\vert t\right\vert ^{-\delta }$ and $_{-\infty }D_{t}^{\alpha
}\left\vert t\right\vert ^{-\delta }$, as~well as to justify that
\begin{equation*}
_{-\infty }D_{t}^{\alpha }\left\vert t\right\vert ^{-\delta }=\,_{-\infty
}I_{t}^{-\alpha }\left\vert t\right\vert ^{-\delta }.
\end{equation*}
\vsn
This paper is organized as follows. Section~\ref{Section: Preliminaries}\
collects all the definitions of the special functions and polynomials that
appear throughout the paper. Section~\ref{Section: Fractional Integrals}
calculates the fractional integrals $_{0}I_{t}^{\alpha }\left[ t^{\nu
-1}\log t\right] $, $_{0}I_{t}^{\alpha }e^{\lambda t}$, and~$_{-\infty
}I_{t}^{\alpha }\left\vert t\right\vert ^{-\delta }$. Despite the fact that $%
_{0}I_{t}^{\alpha }\left[ t^{\nu -1}\log t\right] $, and~$_{0}I_{t}^{\alpha
}e^{\lambda t}$ are found in the existing literature, it is worth performing these
calculations, as they will be useful in the following section. In Section \ref%
{Section: Fractional Derivatives},\ we calculate $_{0}D_{t}^{\alpha }\left[
t^{\nu -1}\log t\right] $, $_{0}D_{t}^{\alpha }e^{\lambda t}$ and $_{-\infty
}D_{t}^{\alpha }\left\vert t\right\vert^{-\delta }$. Finally, we
collect our conclusions in Section~\ref{Section: Conclusions}.

\section{Preliminaries \label{Section: Preliminaries}}
In this section, we collect all the definitions of the special functions and
polynomials that appear throughout the~paper.
\vsn
{\bf Definition 3.}
For $\mathrm{Re}\,\alpha >0$, the~gamma function is defined as \cite{Lebedev} [Eqn.
1.1.1]%
\begin{equation}
\Gamma \left( z\right) =\int_{0}^{\infty }t^{z-1}e^{-t}dt.  \label{Gamma_def}
\end{equation}
\vsn
{\bf Definition 4.}
The digamma function is defined as \cite{Lebedev} [Eqn. 1.3.1]
\begin{equation}
\psi \left( z\right) =\frac{\Gamma ^{\prime }\left( z\right) }{\Gamma \left(
z\right) }.  \label{digamma_def}
\end{equation}
\vsn
{\bf Definition 5.}
The Pochhammer polynomial is defined as \cite{Atlas} [Eqn. 18:12:1]%
\begin{equation}
\left( x\right) _{n}=x\,\left( x+1\right) \left( x+2\right) \cdots \left(
x+n-1\right) =\frac{\Gamma \left( x+n\right) }{\Gamma \left( x\right) }.
\label{Pochhammer_def}
\end{equation}
\vsn
{\bf Definition 6.}
For $\mathrm{Re\,}a>0$, $\mathrm{Re\,}b>0$, the~beta function is defined as
\cite{Lebedev} [Eqns. 1.5.2\&5]\
\begin{equation}
\mathrm{B}\left( a,b\right) =\int_{0}^{1}t^{a-1}\left( 1-t\right) ^{b-1}dt=%
\frac{\Gamma \left( a\right) \,\Gamma \left( b\right) }{\Gamma \left(
a+b\right) }.  \label{beta_definition}
\end{equation}
\vsn
{\bf  Remark 3.}
\label{Remark: beta}Note that if we are considering $a,b\in
\mathbb{R}
$, then the condition $\mathrm{Re\,}a>0$, $\mathrm{Re\,}b>0$, becomes$%
\mathrm{\,}a,b>0$.
\vsn
{\bf Definition 7.}
The lower incomplete gamma function is defined as \cite{Atlas} [Eqn. 45:3:1]
\begin{equation}
\gamma \left( \alpha ,z\right) =\int_{0}^{z}t^{\alpha -1}e^{-t}dt.
\label{gamma_def}
\end{equation}
\vsn
{\bf Definition 8.}
The two-parameter Mittag--Leffler function is defined as \cite
{NIST} [Eqn. 10.46.3],%
\begin{equation}
E_{\mu ,\nu }\left( z\right) =\sum_{k=0}^{\infty }\frac{z^{k}}{\Gamma \left(
\mu \,k+\nu \right) }.  \label{ML_def}
\end{equation}
\vsn
{\bf Definition 9}
The generalized hypergeometric function is defined as \cite
{NIST} [Eqn. 16.2.1]%
\begin{equation}
_{p}F_{q}\left( \left.
\begin{array}{c}
a_{1},\ldots ,a_{p} \\
b_{1},\ldots ,b_{q}%
\end{array}%
\right\vert z\right) =\sum_{k=0}^{\infty }\frac{\left( a_{1}\right)
_{k}\cdots \left( a_{p}\right) _{k}}{\left( b_{1}\right) _{k}\cdots \left(
b_{q}\right) _{k}}\frac{z^{k}}{k!}.  \label{Hyper_def}
\end{equation}

\section{Fractional~Integrals \label{Section: Fractional Integrals}}
\subsection{Fractional Integral of the Power~Function}
\vsn
{\bf Theorem 1.}
For $\alpha >0$, $\gamma >-1$, and~$t>0$, we have \cite{Ederly} [Eqn. 13.1.(7)]
\begin{equation}
\,_{0}I_{t}^{\alpha }t^{\gamma }=\frac{\Gamma \left( 1+\gamma \right) }{%
\Gamma \left( 1+\gamma +\alpha \right) }t^{\gamma +\alpha }.
\label{Theorem_0}
\end{equation}
\vsn
{\bf Proof.}
Apply definition (\ref{Definition_RL}), and~perform the change in variables $%
t-\tau =st$, with~$t>0$ to~obtain%
\begin{eqnarray*}
_{0}I_{t}^{\alpha }t^{\gamma } &=&\frac{1}{\Gamma \left( \alpha \right) }%
\int_{0}^{t}\left( t-\tau \right) ^{\alpha -1}\tau ^{\gamma }\,d\tau \\
&=&\frac{t^{\alpha +\gamma }}{\Gamma \left( \alpha \right) }\underset{%
\mathrm{B}\left( \alpha ,1+\gamma \right) }{\underbrace{\int_{0}^{1}s^{%
\alpha -1}\left( 1-s\right) ^{\gamma }\,ds }}.
\end{eqnarray*}%
Finally, apply (\ref{beta_definition})\ to obtain (\ref{Theorem_0}) as~we
wanted to prove. Note that we consider $\alpha ,\gamma \in
\mathbb{R}
$; thus, according to Remark \ref{Remark: beta}, we have $\alpha >0$, $\gamma
>-1$.

As mentioned in the Introduction, the~Weyl fractional integral for the power
function given in (\ref{Mainardi_power})\ does not seem to hold true.
Therefore, next, we are going to calculate $_{-\infty }I_{t}^{\alpha
}\left\vert t\right\vert ^{-\delta }$. For~this purpose, let us first prove
the following~lemma.
\vsn
{\bf Lemma 1.}
For $\alpha <-\beta -1<0$, and~$t>0$, the~following integral formula holds
true:%
\begin{eqnarray}
\mathcal{I}_{\alpha ,\beta }\left( t\right) &=&\int_{-\infty }^{0}\left(
t-\tau \right) ^{\alpha }\left\vert \tau \right\vert ^{\beta }d\tau
\label{Lemma_1} \\
&=&t^{\alpha +\beta +1}\frac{\Gamma \left( -1-\alpha -\beta \right) \,\Gamma
\left( \beta +1\right) }{\Gamma \left( -\alpha \right) }.  \notag
\end{eqnarray}
\vsn
{\bf Proof.}
Perform the change in variables $u=-\tau $ to~obtain%
\begin{eqnarray*}
\mathcal{I}_{\alpha ,\beta }\left( t\right) &=&\int_{-\infty }^{0}\left(
t-\tau \right) ^{\alpha }\left( -\tau \right) ^{\beta }d\tau \\
&=&\int_{0}^{\infty }\left( t+u\right) ^{\alpha }u^{\beta }du.
\end{eqnarray*}%
Perform the change in variables $u=st/(1-s)$ with $t>0$, and~use the
definition of the beta function (\ref{beta_definition}) to obtain%
\begin{eqnarray*}
\mathcal{I}_{\alpha ,\beta }\left( t\right) &=&\int_{0}^{1}\left( t+\frac{st%
}{1-s}\right) ^{\alpha }\left( \frac{st}{1-s}\right) ^{\beta }t\frac{ds}{%
\left( 1-s\right) ^{2}} \\
&=&t^{\alpha +\beta +1}\underset{\mathrm{B}\left( -\alpha -\beta -1,\beta
+1\right) }{\underbrace{\int_{0}^{1}\left( 1-s\right) ^{-\alpha -\beta
-2}s^{\beta }ds}}.
\end{eqnarray*}%
Applying again (\ref{beta_definition}), we finally arrive at (\ref{Lemma_1}%
), as~we wanted to prove. According to Remark~\ref{Remark: beta}, we have $%
\alpha <-\beta -1<0$.
\vsn
{\bf Theorem 2,}
For $0<\alpha <\delta <1$, and~$t>0$%
\begin{equation}
_{-\infty }I_{t}^{\alpha }\left\vert t\right\vert ^{-\delta }=\frac{\Gamma
\left( \delta -\alpha \right) }{\Gamma \left( \delta \right) }\frac{\cos
\left( \frac{\pi \delta }{2}-\pi \alpha \right) }{\cos \left( \frac{\pi
\delta }{2}\right) }t^{\alpha -\delta }.  \label{Theorem_1}
\end{equation}
\vsn
{\bf Proof}
Note that%
\begin{equation*}
_{-\infty }I_{t}^{\alpha }\left\vert t\right\vert ^{-\delta }=\,_{-\infty
}I_{0}^{\alpha }\left\vert t\right\vert ^{-\delta }+\,_{0}I_{t}^{\alpha
}\left\vert t\right\vert ^{-\delta }.
\end{equation*}%
Also, from~(\ref{Definition_RL})\ and (\ref{Theorem_0}), for~$\alpha >0$, $%
\delta <1$ and $t>0$, we have
\begin{equation}
\,_{0}I_{t}^{\alpha }\left\vert t\right\vert ^{-\delta }=\,_{0}I_{t}^{\alpha
}\,t^{-\delta }=\frac{\Gamma \left( 1-\delta \right) }{\Gamma \left(
1-\delta +\alpha \right) }t^{\alpha -\delta }.  \label{Part_1}
\end{equation}%
Moreover, from~definition (\ref{Definition_RL})\ and lemma (\ref{Lemma_1}),
we have for $\alpha <\delta <1$, and~$t>0$
\begin{equation}
_{-\infty }I_{0}^{\alpha }\left\vert t\right\vert ^{-\delta }=\frac{1}{%
\Gamma \left( \alpha \right) }\underset{\mathcal{I}_{\alpha -1,-\delta
}\left( t\right) }{\underbrace{\int_{-\infty }^{0}\!\!\left( t-\tau \right)
^{\alpha -1}\left\vert \tau \right\vert ^{-\delta }d\tau}}=t^{\alpha -\delta }%
\frac{\Gamma \left( \delta -\alpha \right) \,\Gamma \left( 1-\delta \right)
}{\Gamma \left( 1-\alpha \right) },  \label{Part_2}
\end{equation}%
and thus, according to (\ref{Part_1})\ and (\ref{Part_2}), we obtain
\begin{eqnarray*}
&&_{-\infty }I_{t}^{\alpha }\left\vert t\right\vert ^{-\delta } \\
&=&t^{\alpha -\delta }\,\Gamma \left( 1-\delta \right) \left[ \frac{\Gamma
\left( \delta -\alpha \right) }{\Gamma \left( \alpha \right) \,\Gamma \left(
1-\alpha \right) }+\frac{1}{\Gamma \left( 1+\alpha -\delta \right) }\right]
\\
&=&t^{\alpha -\delta }\,\frac{\Gamma \left( \delta -\alpha \right) \,\Gamma
\left( 1-\delta \right) \,\Gamma \left( \delta \right) }{\Gamma \left(
\delta \right) } \\
&&\left[ \frac{1}{\Gamma \left( \alpha \right) \,\Gamma \left( 1-\alpha
\right) }+\frac{1}{\Gamma \left( 1+\alpha -\delta \right) \,\Gamma \left(
\delta -\alpha \right) }\right] .
\end{eqnarray*}%
Now, apply the property \cite{Lebedev} [Eqn. 1.2.2]\
\begin{equation}
\Gamma \left( z\right) \,\Gamma \left( 1-z\right) =\frac{\pi }{\sin \left(
\pi z\right) }  \label{Reflection_gamma}
\end{equation}%
to obtain
\begin{equation*}
_{-\infty }I_{t}^{\alpha }\left\vert t\right\vert ^{-\delta }=t^{\alpha
-\delta }\,\frac{\Gamma \left( \delta -\alpha \right) }{\Gamma \left( \delta
\right) \sin \left( \pi \delta \right) }\left[ \sin \left( \pi \alpha
\right) +\sin \left( \pi \left( \delta -\alpha \right) \right) \right] .
\end{equation*}%
Finally, apply the property \cite{Spiegel} [Eqn. 5.61]
\begin{equation}
\sin A+\sin B=2\sin \left( \frac{A+B}{2}\right) \cos \left( \frac{A-B}{2}%
\right) ,  \label{Trig_property}
\end{equation}%
to arrive at (\ref{Theorem_1}) as~we wanted to prove.
\vsn
\subsection{Fractional Integral of the Exponential~Function}
\vsn
{\bf Lemma 2.}
The following identity holds true:%
\begin{equation}
\gamma \left( \alpha ,z\right) =z^{\alpha }\,\Gamma \left( \alpha \right)
\,e^{-z}\,E_{1,1+\alpha }\left( z\right) .  \label{gamma_ML}
\end{equation}
\vsn
{\bf Proof}
Consider the expansion \cite{Atlas} [Eqn. 45:6:2]%
\begin{equation*}
e^{z}\,\gamma \left( \alpha ,z\right) =\frac{z^{\alpha }}{\alpha }%
\sum_{k=0}^{\infty }\frac{z^{k}}{\left( \alpha +1\right) _{k}}.
\end{equation*}%
Taking into account the factorial property of the gamma function \cite{Lebedev} [Eqn.
1.2.1], i.e.,
\begin{equation}
\Gamma \left( z+1\right) =z\,\Gamma \left( z\right) ,
\label{Gamma_factorial}
\end{equation}%
and the definition of the Pochhammer symbol (\ref{Pochhammer_def}), we have%
\begin{equation*}
e^{z}\,\gamma \left( \alpha ,z\right) =z^{\alpha }\,\Gamma \left( \alpha
\right) \sum_{k=0}^{\infty }\frac{z^{k}}{\Gamma \left( \alpha +1+k\right) }.
\end{equation*}%
Finally, apply the definition of the Mittag--Leffler function (\ref{ML_def})\
to complete the proof.
\vsn
{\bf Theorem 3.}
For $t>0$, and~$\alpha >0$, the~following fractional integral holds true
\cite{Magin} [Table IV.1]:%
\begin{equation}
_{0}I_{t}^{\alpha }e^{\lambda t}=t^{\alpha }E_{1,1+\alpha }\left( \lambda
t\right) .  \label{FInt_exp}
\end{equation}
\vsn
{\bf Proof.}
According to the definition of the Riemann--Liouville fractional integral (%
\ref{Definition_RL}), we have%
\begin{eqnarray*}
_{0}I_{t}^{\alpha }e^{\lambda t} &=&\frac{1}{\Gamma \left( \alpha \right) }%
\int_{0}^{t}\left( t-\tau \right) ^{\alpha -1}e^{\lambda \tau }d\tau \\
&=&\frac{\lambda ^{-\alpha }}{\Gamma \left( \alpha \right) }%
\int_{0}^{t}\left( \lambda t-\lambda \tau \right) ^{\alpha -1}e^{\lambda
\tau }\lambda d\tau ,
\end{eqnarray*}%
thus performing the change in variables $u=\lambda \left( t-\tau \right) $,
and taking into account the definition of the lower incomplete gamma
function (\ref{gamma_def}), we obtain%
\begin{eqnarray*}
_{0}I_{t}^{\alpha }e^{\lambda t} &=&\frac{\lambda ^{-\alpha }e^{\lambda t}}{%
\Gamma \left( \alpha \right) }\int_{0}^{\lambda t}u^{\alpha -1}e^{-u}du \\
&=&\frac{\lambda ^{-\alpha }e^{\lambda t}}{\Gamma \left( \alpha \right) }%
\,\gamma \left( \alpha ,\lambda t\right) .
\end{eqnarray*}%
Finally, apply (\ref{gamma_ML}) to complete the proof.
\vsn
\subsection{Fractional Integral Formula of the Logarithmic~Function}
\vsn
{\bf Lemma 3.}
The following integral formula holds true:%
\begin{eqnarray}
&&\int_{0}^{1}t^{a-1}\left( 1-t\right) ^{b-1}\log t\,dt=\frac{\Gamma \left(
a\right) \,\Gamma \left( b\right) }{\Gamma \left( a+b\right) }\left[ \psi
\left( a\right) -\psi \left( a+b\right) \right] ,  \label{D_beta} \\
&&\mathrm{Re\,}a>0,\,\mathrm{Re\,}b>0.  \notag
\end{eqnarray}
\vsn
{\bf Proof.}
Perform the derivative with respect to the first parameter in the beta
function definition (\ref{beta_definition}),
\begin{eqnarray}
&&\frac{\partial \,}{\partial a}\mathrm{B}\left( a,b\right)  \label{D_beta_2}
\\
&=&\int_{0}^{1}t^{a-1}\left( 1-t\right) ^{b-1}\log t\,dt=\,\Gamma \left(
b\right) \,\frac{\Gamma ^{\prime }\left( a\right) \,\Gamma \left( a+b\right)
-\Gamma \left( a\right) \,\Gamma ^{\prime }\left( a+b\right) }{\Gamma
^{2}\left( a+b\right) },  \notag
\end{eqnarray}%
and take into account the definition of the digamma function (\ref%
{digamma_def}) to complete the proof.
\vsn
{\bf Theorem 4.}
The following fractional integral holds true \cite{Ederly} [Eqn. 13.1.(24)]:\
\begin{eqnarray}
&&_{0}I_{t}^{\alpha }\left[ t^{\nu -1}\log t\right] =\frac{t^{\alpha +\nu
-1}\,\,\Gamma \left( \nu \right) }{\Gamma \left( \alpha +\nu \right) }\left[
\log t+\psi \left( \nu \right) -\psi \left( \alpha +\nu \right) \right] ,
\label{I_log(t)} \\
&&\mathrm{Re\,}\alpha >0,\,\mathrm{Re\,}\nu >0.  \notag
\end{eqnarray}
\vsn
{\bf Proof.}
Apply the definition (\ref{Definition_RL})\ and perform the change in
variables $\tau =t\,u$ with $t>0$ to~obtain%
\begin{eqnarray}
&&_{0}I_{t}^{\alpha }\left[ t^{\nu -1}\log t\right]  \notag \\
&=&\frac{1}{\Gamma \left( \alpha \right) }\int_{0}^{t}\left( t-\tau \right)
^{\alpha -1}\tau ^{\nu -1}\log \tau \,d\tau  \notag \\
&=&\frac{t^{\alpha +\nu -1}}{\Gamma \left( \alpha \right) }\left\{ \log
t\,\int_{0}^{1}\left( 1-u\right) ^{\alpha -1}u^{\nu -1}du+\int_{0}^{1}\left(
1-u\right) ^{\alpha -1}u^{\nu -1}\log u\,du\right\} .  \label{Int_1_2}
\end{eqnarray}%
Use the definition of the beta function (\ref{beta_definition})\ to
calculate the first integral in (\ref{Int_1_2})\ as
\begin{eqnarray}
&&\int_{0}^{1}\left( 1-u\right) ^{\alpha -1}u^{\nu -1}du=\frac{\Gamma \left(
\nu \right) \,\Gamma \left( \alpha \right) }{\Gamma \left( \alpha +\nu
\right) },  \label{Int_1} \\
&&\mathrm{Re\,}\alpha >0,\,\mathrm{Re\,}\nu >0.  \notag
\end{eqnarray}%
The second integral in (\ref{Int_1_2})\ is given in lemma (\ref{D_beta}).
Therefore, substituting (\ref{D_beta})\ and (\ref{Int_1})\ in (\ref{Int_1_2}%
), we arrive at (\ref{I_log(t)}) as~we wanted to prove.
\section{Fractional~Derivatives \label{Section: Fractional Derivatives}}
\subsection{Fractional Derivative of the Power~Function}
\vsn
{\bf Lemma 4.}
The following $n$-th derivative formula holds:%
\begin{equation}
\frac{d^{n}}{dt^{n}}t^{a}=\frac{\Gamma \left( a+1\right) }{\Gamma \left(
a-n+1\right) }t^{a-n}.  \label{Dn_t^a}
\end{equation}
\vsn
{Proof.}
According to the definition of the Pochhamer symbol (\ref{Pochhammer_def}),
we have%
\begin{eqnarray*}
\frac{d^{n}}{dt^{n}}t^{a} &=&a\left( a-1\right) \cdots \left( a-n+1\right)
t^{a-n} \\
&=&\left( a-n+1\right) _{n}\,t^{a-n} \\
&=&\frac{\Gamma \left( a+1\right) }{\Gamma \left( a-n+1\right) }t^{a-n}.
\end{eqnarray*}
\vsn
{\bf Theorem 5.}
For $\alpha >0$, $\gamma >-1$, and~$t>0$ \cite{Magin} [Table IV.1], \cite{Polubny} [%
Appendix], the~following fractional derivative holds true:%
\begin{equation}
_{0}D_{t}^{\alpha }t^{\gamma }=\frac{\Gamma \left( 1+\gamma \right) }{\Gamma
\left( 1+\gamma -\alpha \right) }t^{\gamma -\alpha }.  \label{D_RL_t^a}
\end{equation}
\vsn
{\bf Proof.}
According to definition (\ref{Derivative_RL_def}), and~the results given in (%
\ref{Theorem_0})\ and (\ref{Dn_t^a}),
\begin{eqnarray*}
_{0}D_{t}^{\alpha }t^{\gamma } &=&\left( D_{t}^{m}\circ
\,_{0}I_{t}^{m-\alpha }\right) t^{\gamma } \\
&=&\frac{\Gamma \left( 1+\gamma \right) }{\Gamma \left( 1+\gamma +m-\alpha
\right) }\frac{d^{m}}{dt^{m}}t^{\gamma +m-\alpha } \\
&=&\frac{\Gamma \left( 1+\gamma \right) }{\Gamma \left( 1+\gamma +m-\alpha
\right) }\frac{\Gamma \left( \gamma +m-\alpha +1\right) }{\Gamma \left(
\gamma -\alpha +1\right) }t^{\gamma +m-\alpha -m},
\end{eqnarray*}%
as we wanted to prove.
\vsn

Now, we calculate the Weyl fractional derivative corresponding to the Weyl
fractional integral calculated in (\ref{Theorem_1}).
\vsn
{\bf Theorem 6.}
For $0<\delta <1$, $\alpha >0$, and~$t>0$,
\begin{equation}
_{-\infty }D_{t}^{\alpha }\left\vert t\right\vert ^{-\delta }=\frac{\Gamma
\left( \delta +\alpha \right) }{\Gamma \left( \delta \right) }\frac{\cos
\left( \frac{\pi \delta }{2}+\pi \alpha \right) }{\cos \left( \frac{\pi
\delta }{2}\right) }t^{-\alpha -\delta }.  \label{Theorem_2}
\end{equation}
\vsn
{\bf Proof.}
According to the definition of the Weyl fractional integral (\ref%
{Derivative_Weyl_def}), we have \
\begin{equation}
_{-\infty }D_{t}^{\alpha }\left\vert t\right\vert ^{-\delta }=\,_{-\infty
}D_{0}^{\alpha }\left\vert t\right\vert ^{-\delta }+\,_{0}D_{t}^{\alpha
}\left\vert t\right\vert ^{-\delta }.  \label{Sum_D}
\end{equation}%
On the one hand, apply definition (\ref{Derivative_RL_def}), taking into
account that $\alpha \notin
\mathbb{N}
$ and $m-\alpha >0$,
\begin{equation*}
\,_{-\infty }D_{0}^{\alpha }\left\vert t\right\vert ^{-\delta }=\frac{1}{%
\Gamma \left( m-\alpha \right) }\frac{d^{m}}{dt^{m}}\underset{\mathcal{I}%
_{m-1-\alpha ,-\delta }\left( t\right) }{\underbrace{\int_{-\infty }^{0}%
\frac{\left\vert \tau \right\vert ^{-\delta }}{\left( t-\tau \right)
^{\alpha +1-m}}d\tau }}.
\end{equation*}%
Now, apply lemmas (\ref{Lemma_1}) and (\ref{Dn_t^a}) to~obtain for $%
m-\alpha <\delta <1$%
\begin{eqnarray}
_{-\infty }D_{0}^{\alpha }\left\vert t\right\vert ^{-\delta } &=&\frac{%
\Gamma \left( \alpha +\delta -m\right) \,\Gamma \left( 1-\delta \right) }{%
\Gamma \left( m-\alpha \right) \,\Gamma \left( 1+\alpha -m\right) }\frac{%
d^{m}}{dt^{m}}t^{m-\alpha -\delta }  \notag \\
&=&\frac{\Gamma \left( \alpha +\delta -m\right) \,\Gamma \left( 1-\delta
\right) }{\Gamma \left( m-\alpha \right) \,\Gamma \left( 1+\alpha -m\right) }%
\frac{\Gamma \left( m-\alpha -\delta +1\right) }{\Gamma \left( -\alpha
-\delta +1\right) }t^{-\alpha -\delta }.  \label{Part_A}
\end{eqnarray}%
On the other hand, taking into account that $t>0$ and $\alpha \notin
\mathbb{N}
$, apply (\ref{D_RL_t^a})\ to obtain for $\delta <1$,%
\begin{equation}
_{0}D_{t}^{\alpha }\left\vert t\right\vert ^{-\delta }=\,_{0}D_{t}^{\alpha
}\,t^{-\delta }=\frac{\Gamma \left( 1-\delta \right) }{\Gamma \left(
1-\delta -\alpha \right) }t^{-\delta -\alpha }.  \label{Part_B}
\end{equation}%
Substitute results (\ref{Part_A})\ and (\ref{Part_B})\ in (\ref{Sum_D}) to~arrive at%
\begin{equation*}
_{-\infty }D_{t}^{\alpha }\left\vert t\right\vert ^{-\delta }=\frac{%
t^{-\alpha -\delta }\,\Gamma \left( 1-\delta \right) }{\Gamma \left(
1-\alpha -\delta \right) }\left[ \frac{\Gamma \left( \alpha +\delta
-m\right) \,\Gamma \left( 1-\alpha -\delta +m\right) }{\Gamma \left(
m-\alpha \right) \,\Gamma \left( 1+\alpha -m\right) }+1\right] .
\end{equation*}%
Now, apply (\ref{Reflection_gamma})%
\begin{eqnarray*}
&&_{-\infty }D_{t}^{\alpha }\left\vert t\right\vert ^{-\delta } \\
&=&t^{-\alpha -\delta }\frac{\,\Gamma \left( \alpha +\delta \right) }{\Gamma
\left( \delta \right) }\frac{\,\Gamma \left( 1-\delta \right) \,\Gamma
\left( \delta \right) }{\Gamma \left( 1-\alpha -\delta \right) \,\Gamma
\left( \alpha +\delta \right) }\left[ \frac{\sin \pi \left( m-\alpha \right)
}{\sin \pi \left( \alpha +\delta -m\right) }+1\right] \\
&=&t^{-\alpha -\delta }\frac{\,\Gamma \left( \alpha +\delta \right) }{\Gamma
\left( \delta \right) }\frac{\sin \pi \left( \alpha +\delta \right) }{\sin
\pi \delta }\left[ 1-\frac{\sin \pi \alpha }{\sin \pi \left( \alpha +\delta
\right) }\right] \\
&=&t^{-\alpha -\delta }\frac{\,\Gamma \left( \alpha +\delta \right) }{\Gamma
\left( \delta \right) }\left[ \frac{\sin \pi \left( \alpha +\delta \right)
-\sin \pi \alpha }{\sin \pi \delta }\right] .
\end{eqnarray*}%
Finally, apply (\ref{Trig_property})\ and simplify to arrive at (\ref%
{Theorem_2}) as~we wanted to prove.
\vsn
\subsection{Fractional Derivative of the Exponential~Function}
\vsn
{\bf Theorem 7.}
For $\alpha >0$, $\alpha \notin
\mathbb{N}
$, and~$t>0$%
\begin{equation}
_{0}D_{t}^{\alpha }e^{\lambda t}=t^{-\alpha }E_{1,1-\alpha }\left( \lambda
t\right) .  \label{D_exp}
\end{equation}
\vsn
{\bf Proof.}
Apply the definition (\ref{Derivative_RL_def})\ and expand the exponential
fraction in its Maclaurin series to obtain%
\begin{eqnarray}
_{0}D_{t}^{\alpha }e^{\lambda t} &=&\frac{1}{\Gamma \left( m-\alpha \right) }%
\frac{d^{m}}{dt^{m}}\int_{0}^{t}\left( t-\tau \right) ^{m-\alpha
-1}e^{\lambda \tau }d\tau  \notag \\
&=&\frac{1}{\Gamma \left( m-\alpha \right) }\frac{d^{m}}{dt^{m}}%
\int_{0}^{t}\left( t-\tau \right) ^{m-\alpha -1}\sum_{k=0}^{\infty }\frac{%
\left( \lambda \tau \right) ^{k}}{k!}d\tau  \notag \\
&=&\frac{1}{\Gamma \left( m-\alpha \right) }\frac{d^{m}}{dt^{m}}\left(
\sum_{k=0}^{\infty }\frac{\lambda ^{k}}{k!}\int_{0}^{t}\frac{\tau ^{k}}{%
\left( t-\tau \right) ^{\alpha -m+1}}d\tau \right) .  \label{Int_beta}
\end{eqnarray}%
Perform the change in variables $\tau =u\,t$ with $t>0$ in (\ref{Int_beta}%
), and~apply the definition of the beta function (\ref{beta_definition}). Thus,
for $m-\alpha >0$, we have
\begin{eqnarray*}
_{0}D_{t}^{\alpha }e^{\lambda t} &=&\frac{1}{\Gamma \left( m-\alpha \right) }%
\frac{d^{m}}{dt^{m}}\left( t^{m-\alpha }\sum_{k=0}^{\infty }\frac{\left(
\lambda t\right) ^{k}}{k!}\int_{0}^{1}u^{k}\,\left( 1-u\right) ^{m-\alpha
-1}du\right) \\
&=&\frac{1}{\Gamma \left( m-\alpha \right) }\frac{d^{m}}{dt^{m}}\left(
t^{m-\alpha }\sum_{k=0}^{\infty }\frac{\left( \lambda t\right) ^{k}}{k!}\,%
\mathrm{B}\left( k+1,m-\alpha \right) \right) \\
&=&\frac{d^{m}}{dt^{m}}\sum_{k=0}^{\infty }\frac{\lambda ^{k}\,t^{m-\alpha
+k}}{\Gamma \left( m-\alpha +k+1\right) }.
\end{eqnarray*}%
Now, apply the differentiation formula (\ref{Dn_t^a}) to~arrive at%
\begin{equation*}
_{0}D_{t}^{\alpha }e^{\lambda t}=t^{-\alpha }\sum_{k=0}^{\infty }\frac{%
\left( \lambda t\right) ^{k}}{\Gamma \left( k+1-\alpha \right) }.
\end{equation*}%
Finally, apply the definition of the Mittag--Leffler function (\ref{ML_def})\
to complete the proof.
\vsn
\subsection{Fractional Derivative of the Logarithm~Function}
\vsn
{\bf Lemma 5.}
The following $n$-th derivative formula holds true:%
\begin{equation}
\frac{d^{n}}{dt^{n}}\left( t^{\beta }\log t\right) =\frac{\Gamma \left(
\beta +1\right) }{\Gamma \left( \beta -n+1\right) }t^{\beta -n}\left[ \log
t+\psi \left( \beta +1\right) -\psi \left( \beta -n+1\right) \right] .
\label{Dn_log*t^a}
\end{equation}
\vsn
{\bf Proof.}
According to Leibniz's differentiation formula \cite{NIST} [Eqn. 1.4.12],
\begin{equation}
\frac{d^{n}}{dx^{n}}\left[ f\left( x\right) g\left( x\right) \right]
=\sum_{k=0}^{n}\binom{n}{k}\frac{d^{n-k}}{dx^{n-k}}f\left( x\right) \,\frac{%
d^{k}}{dx^{k}}g\left( x\right) ,  \label{Leibniz_formula}
\end{equation}%
and applying the $n$-th derivative formula given in (\ref{Dn_t^a}), as~well
as the following one (which can be easily proved by induction)%
\begin{equation*}
\frac{d^{n}}{dt^{n}}\log t=\frac{\left( -1\right) ^{n-1}\left( n-1\right) !}{%
t^{n}},\quad n=1,2,\ldots
\end{equation*}%
after simplification, we arrive at
\begin{eqnarray}
&&\frac{d^{n}}{dt^{n}}\left( t^{\beta }\log t\right)  \notag \\
&=&\sum_{k=0}^{n}\binom{n}{k}\frac{d^{n-k}}{dt^{n-k}}t^{\beta }\,\frac{d^{k}%
}{dt^{k}}\log t  \notag \\
&=&\log t\,\frac{d^{n}}{dt^{n}}t^{\beta }\,+\sum_{k=1}^{n}\binom{n}{k}\frac{%
d^{n-k}}{dt^{n-k}}t^{\beta }\,\frac{d^{k}}{dt^{k}}\log t  \notag \\
&=&\Gamma \left( \beta +1\right) \,t^{\beta -n}\left[ \frac{\log t}{\Gamma
\left( \beta -n+1\right) }+n!\sum_{k=1}^{n}\frac{\left( -1\right) ^{k-1}}{%
k\,\left( n-k\right)! \, \Gamma \left( \beta -n+k+1\right) }\right] .
\label{S_def}
\end{eqnarray}%
In order to calculate the finite sum given in (\ref{S_def}), consider the
following function $f\left( z\right) $ which can be recast as a hypergeometric function \cite{Andrews} [Sect. 2.1]
\begin{eqnarray}
f\left( z\right) &=&\sum_{k=1}^{n}\frac{z^{k-1}}{\left( n-k\right) !\,\Gamma
\left( \beta -n+k+1\right) }  \label{f(z)_1} \\
&=&\frac{1}{\Gamma \left( n\right) \Gamma \left( \beta -n+2\right) }%
\,_{2}F_{1}\left( \left.
\begin{array}{c}
1,1-n \\
\beta -n+2%
\end{array}%
\right\vert -z\right) .  \label{f(z)_2}
\end{eqnarray}%
On the one hand, integrating term by term in (\ref{f(z)_1})%
\begin{equation}
g\left( z\right) =\int_{0}^{z}f\left( t\right) \,dt=\sum_{k=1}^{n}\frac{z^{k}%
}{k\,\left( n-k\right) !\,\Gamma \left( \beta -n+k+1\right) }.
\label{g(z)_1}
\end{equation}%
On the other hand, applying the integration formula given in \cite[Eqn. 2.2.3%
]{Andrews}%
\begin{eqnarray*}
&&_{p+1}F_{q+1}\left( \left.
\begin{array}{c}
a_{1},\ldots ,a_{p},a_{p+1} \\
b_{1},\ldots ,b_{q},b_{q+1}%
\end{array}%
\right\vert x\right) \\
&=&\frac{\Gamma \left( b_{q+1}\right) \,x^{1-b_{q+1}}}{\Gamma \left(
a_{p+1}\right) \,\Gamma \left( b_{q+1}-a_{p+1}\right) } \\
&&\int_{0}^{x}t^{a_{p+1}-1}\left( x-t\right)
^{b_{q+1}-a_{p+1}-1}\,_{p}F_{q}\left( \left.
\begin{array}{c}
a_{1},\ldots ,a_{p} \\
b_{1},\ldots ,b_{q}%
\end{array}%
\right\vert t\right) \,dt,
\end{eqnarray*}%
taking $a_{1}=1$, $a_{2}=1-n$, $b_{1}=\beta -n+2$, we obtain%
\begin{eqnarray}
&&g\left( z\right)  \notag \\
&=&\int_{0}^{z}f\left( t\right) \,dt=\frac{1}{\Gamma \left( n\right) \Gamma
\left( \beta -n+2\right) }\int_{0}^{z}\,_{2}F_{1}\left( \left.
\begin{array}{c}
1,1-n \\
\beta -n+2%
\end{array}%
\right\vert -t\right) \,dt  \notag \\
&=&\frac{z}{\Gamma \left( n\right) \Gamma \left( \beta -n+2\right) }%
\,_{3}F_{2}\left( \left.
\begin{array}{c}
1,1,1-n \\
2,\beta -n+2%
\end{array}%
\right\vert -z\right)  \label{g(z)_2}
\end{eqnarray}%
Therefore, from~(\ref{g(z)_1})\ and (\ref{g(z)_2}),
\begin{eqnarray*}
-g\left( -1\right) &=&\sum_{k=1}^{n}\frac{\left( -1\right) ^{k-1}}{k\,\left(
n-k\right) !\,\Gamma \left( \beta -n+k+1\right) } \\
&=&\frac{1}{\Gamma \left( n\right) \Gamma \left( \beta -n+2\right) }%
\,_{3}F_{2}\left( \left.
\begin{array}{c}
1,1,1-n \\
2,\beta -n+2%
\end{array}%
\right\vert 1\right) .
\end{eqnarray*}%
Now, consider the reduction formula \cite{Prudnikov3} [Eqn. 7.4.4(40)]
\begin{equation*}
_{3}F_{2}\left( \left.
\begin{array}{c}
1,1,a \\
2,b%
\end{array}%
\right\vert 1\right) =\frac{b-1}{a-1}\left[ \psi \left( b-1\right) -\psi
\left( b-a\right) \right] ,
\end{equation*}%
to arrive at%
\begin{equation}
\sum_{k=1}^{n}\frac{\left( -1\right) ^{k-1}}{k\,\left( n-k\right) !\Gamma
\left( \beta -n+k+1\right) }=\frac{\psi \left( \beta +1\right) -\psi \left(
\beta -n+1\right) }{n!\,\Gamma \left( \beta -n+1\right) }.
\label{S_resultado}
\end{equation}%
Finally, substitute (\ref{S_resultado})\ in (\ref{S_def})\ to arrive at the
desired result.
\vsn
{\bf Theorem 8.}
For $\alpha >0$, $\alpha \notin
\mathbb{N}
$, $\mathrm{Re}\,\nu >0$, and~$t>0$, the~following fractional derivative
holds true:%
\begin{equation}
_{0}D_{t}^{\alpha }\left[ t^{\nu -1}\log t\right] =\frac{t^{\nu -\alpha
-1}\,\,\Gamma \left( \nu \right) }{\Gamma \left( \nu -\alpha \right) }\left[
\log t+\psi \left( \nu \right) -\psi \left( \nu -\alpha \right) \right] .
\label{D_log(t)}
\end{equation}
\vsn
{\bf Proof.}
Apply the definition (\ref{Derivative_RL_def})\ and perform the change in
variables $\tau =t\,u$ with $t>0$, to~obtain for $\alpha >0$,%
\begin{eqnarray}
&&_{0}D_{t}^{\alpha }\left[ t^{\nu -1}\log t\right]  \notag \\
&=&\frac{1}{\Gamma \left( m-\alpha \right) }\frac{d^{m}}{dt^{m}}%
\int_{0}^{t}\left( t-\tau \right) ^{m-\alpha -1}\tau ^{\nu -1}\log \tau
\,d\tau  \notag \\
&=&\frac{1}{\Gamma \left( m-\alpha \right) }\frac{d^{m}}{dt^{m}}\left\{
t^{m-\alpha +\nu -1}\left[ \log t\int_{0}^{1}\left( 1-u\right) ^{m-\alpha
-1}u^{\nu -1}du\right. \right.  \notag \\
&&\quad \qquad \qquad +\left. \left. \int_{0}^{1}\left( 1-u\right)
^{m-\alpha -1}u^{\nu -1}\log u\,du\right] \right\}.  \label{Int_log}
\end{eqnarray}%
The first integral in (\ref{Int_log})\ is just a beta function (\ref%
{beta_definition}); thus, for $\mathrm{Re}\,\nu >0$, we have
\begin{equation*}
\int_{0}^{1}\left( 1-u\right) ^{m-\alpha -1}u^{\nu -1}du=\frac{\Gamma \left(
\nu \right) \,\Gamma \left( m-\alpha \right) }{\Gamma \left( m-\alpha +\nu
\right) },
\end{equation*}%
and for the second one, we can apply (\ref{D_beta}); thus, for $\mathrm{Re}%
\,\nu >0$, we have
\begin{equation*}
\int_{0}^{1}\left( 1-u\right) ^{m-\alpha -1}u^{\nu -1}\log u\,du=\frac{%
\Gamma \left( \nu \right) \,\Gamma \left( m-\alpha \right) }{\Gamma \left(
m-\alpha +\nu \right) }\left[ \psi \left( \nu \right) -\psi \left( m-\alpha
+\nu \right) \right] ,
\end{equation*}%
thereby%
\begin{eqnarray}
&&_{0}D_{t}^{\alpha }\left[ t^{\nu -1}\log t\right]  \notag \\
&=&\frac{\Gamma \left( \nu \right) \,}{\Gamma \left( m-\alpha +\nu \right) }
\notag \\
&&\left\{ \frac{d^{m}}{dt^{m}}\left( t^{m-\alpha +\nu -1}\log t\right) +%
\left[ \psi \left( \nu \right) -\psi \left( m-\alpha +\nu \right) \right]
\frac{d^{m}}{dt^{m}}t^{m-\alpha +\nu -1}\right\} .  \label{D_1&2}
\end{eqnarray}%
Apply (\ref{Dn_log*t^a})\ to obtain%
\begin{eqnarray}
&&\frac{d^{m}}{dt^{m}}\left( t^{m-\alpha +\nu -1}\log t\right)  \notag \\
&=&\frac{\Gamma \left( m-\alpha +\nu \right) }{\Gamma \left( \nu -\alpha
\right) }t^{-\alpha +\nu -1}\left[ \log t+\psi \left( m-\alpha +\nu \right)
-\psi \left( \nu -\alpha \right) \right] ,  \label{D_1}
\end{eqnarray}%
and apply (\ref{Dn_t^a})\ to obtain %
\begin{equation}
\frac{d^{m}}{dt^{m}}t^{m-\alpha +\nu -1}=\frac{\Gamma \left( m-\alpha +\nu
\right) }{\Gamma \left( \nu -\alpha \right) }t^{-\alpha +\nu -1}.
\label{D_2}
\end{equation}%
Insert (\ref{D_1})\ and (\ref{D_2})\ in (\ref{D_1&2}), and~simplify the
result to complete the proof.
\vsn
\section{Conclusions \label{Section: Conclusions}}
\vsn
On the one hand, according to (\ref{I_log(t)}), (\ref{D_log(t)}), (\ref%
{FInt_exp})\ and (\ref{D_exp}), we have analytically justified~that%
\begin{eqnarray}
_{0}D_{t}^{\alpha }\left[ t^{\nu -1}\log t\right] &=&\,_{0}I_{t}^{-\alpha }%
\left[ t^{\nu -1}\log t\right]  \label{Log_a->-a} \\
&=&\frac{t^{\nu -\alpha -1}\,\,\Gamma \left( \nu \right) }{\Gamma \left( \nu
-\alpha \right) }\left[ \log t+\psi \left( \nu \right) -\psi \left( \nu
-\alpha \right) \right] ,  \notag \\
_{0}D_{t}^{\alpha }e^{\lambda t} &=&\,_{0}I_{t}^{-\alpha }e^{\lambda
t}=t^{-\alpha }E_{1,1-\alpha }\left( \lambda t\right) \,,  \label{Exp_a->-a}
\end{eqnarray}%
applying the corresponding definitions of the Riemann--Liouville fractional
integral (\ref{Definition_RL}) and~the Riemann--Liouville fractional
derivative (\ref{Derivative_RL_def}). Note that\ the fractional derivatives
calculated in (\ref{Log_a->-a})\ and (\ref{Exp_a->-a})\ can be obtained from the
corresponding fractional integrals, substituting $\alpha $ by $-\alpha $. However, the~corresponding derivations from the Riemann--Liouville definitions of the
fractional integral and the fractional derivative in order to arrive to this
conclusion are highly~non-trivial.

On the other hand, from~the definitions of the Weyl fractional integral (\ref%
{Definition_Weyl_FI}),\ and the Weyl fractional derivative (\ref%
{Derivative_Weyl_def}), we have calculated the novel formulas (\ref%
{Theorem_1}), and~(\ref{Theorem_2}), i.e.,
\begin{equation*}
_{-\infty }D_{t}^{\alpha }\left\vert t\right\vert ^{-\delta }=\frac{\Gamma
\left( \delta +\alpha \right) }{\Gamma \left( \delta \right) }\frac{\cos
\left( \frac{\pi \delta }{2}+\pi \alpha \right) }{\cos \left( \frac{\pi
\delta }{2}\right) }t^{-\alpha -\delta }=\,_{-\infty }I_{t}^{-\alpha
}\left\vert t\right\vert ^{-\delta }.
\end{equation*}%
Again, we can obtain the derivative formula $_{-\infty }D_{t}^{\alpha
}\left\vert t\right\vert ^{-\delta }$ substituting $\alpha $ by $-\alpha $
in the corresponding formula for $_{-\infty }I_{t}^{\alpha }\left\vert
t\right\vert ^{-\delta }$. Nevertheless, according to the corresponding
derivations, this property is not straightforward as~in the case of (\ref%
{Log_a->-a})\ and (\ref{Exp_a->-a}). In~general, this occurs because the
definition of the Riemann--Liouville fractional derivative (\ref%
{Derivative_RL_def}) involves an $m$-th derivative. Meanwhile, this is not the
case for the definition of the Riemann--Liouville fractional integral (\ref%
{Definition_RL}). It would be interesting to investigate the conditions
under which it is satisfied that%
\begin{equation}
_{a}D_{t}^{\alpha }f\left( t\right) =\,_{a}I_{t}^{-\alpha }f\left( t\right) ,
\label{Property_FD_FI}
\end{equation}%
from the definitions given in (\ref{Definition_RL})\ and (\ref%
{Derivative_RL_def}) since (\ref{Property_FD_FI})\ is implicitly taken for
granted in the most common literature about fractional calculus, to the knowledge of the~authors.



\end{document}